# Three-point Step Size Gradient Method with Relaxed Generalized Armijo Step Size Rule[*]

Sun Qingying   Zhao Xu   Wang Jian

**College of Science, China University of Petroleum (East China), Qingdao, Shandong, 266580**

## Abstract

Based on differences of points and differences of gradients over the most recent three iterations, together with the Taylor's theorem, two forms of the quasi-Newton equations at the recent iteration are constructed. By using the two forms of the quasi-Newton equation and the method of least squares, three-point step size gradient methods for solving unconstrained optimization problem are proposed. It is proved by using the relaxed generalized Armijo step size rule that the new method is of global convergence properties if the gradient function is uniformly continuous. Moreover, it is shown that, when the objective function is pseudo-convex (quasi-convex) function, the new method has strong convergence results. In addition, it is also shown under some suitable assumptions that the new method is of super-linear and linear convergence. Although multi-piont information is used, TBB has the feature of simplicity, low memory requirement and only first order information being used, the new method is very suitable for solving large-scale optimization problems. Numerical experiments are provided and the efficiency, robustness and analysis of TBB are confirmed.

**Key words:** Unconstrained optimization, BB method, Three-point step size gradient method, Relaxed generalized Armijo step size rule，Convergence

## 1. Introduction

Consider the unconstrained optimization problem:

$$\min_{x \in R^n} f(x), \qquad (1.1)$$

where $f(x): R^n \to R^1$ is a real continuously differential function.

This problem can be solved by the quasi-Newton method with simple structure and fast convergence, whose iterative formula is

$$x_{k+1} = x_k + \alpha_k d_k, \qquad (1.2)$$

where $\alpha_k$ is the line search step size which can be determined by certain strategy and the direction of search is

$$d_k = -H_k g_k, \qquad (1.3)$$

with $g_k = \nabla f(x_k)$ being the gradient of the objective function at $x_k$, $H_k$ being a kind of approximation to the Hessian matrix of objective function. $H_k$ is required positive definite and to satisfy the quasi Newton equation:

$$H_k y_{k-1} = s_{k-1}, \qquad (1.4)$$

Supported by Chinese Natural Science Foundation with grant 51974343. Corresponding author: Qingying Sun, Email: sunqingying01@163.com

where $y_{k-1} = g_k - g_{k-1}$ and $s_{k-1} = x_k - x_{k-1}$. However, storage of quasi Newton method is at least $O(n^2)$, making it difficult to solve large-scaled problems.

In 1988 Barzilai and Borwein[1]proposed the two-point step size gradient method (called BB method) with the form:

$$x_{k+1} = x_k - D_k g_k , \tag{1.5}$$

where $D_k = \alpha_k I$ is a scalar matrix. In order to make $D_k$ contain second order information of objective function, the step size in computation $\alpha_k$ should approximately satisfy the following quasi Newton equations in the sense of least square:

$$\min_{\alpha>0} \left\| y_{k-1} - \frac{1}{\alpha} s_{k-1} \right\|, \tag{1.6}$$

and

$$\min_{\alpha>0} \left\| \alpha y_{k-1} - s_{k-1} \right\|, \tag{1.7}$$

where $y_{k-1} = g_k - g_{k-1}$ and $s_{k-1} = x_k - x_{k-1}$ are the same as those in (1.4). Solving (1.6) and (1.7) yields

$$\alpha_k^{BB_1} = \frac{\|s_{k-1}\|^2}{s_{k-1}^T y_{k-1}}. \tag{1.8}$$

and

$$\alpha_k^{BB_2} = \frac{s_{k-1}^T y_{k-1}}{\|y_{k-1}\|^2}. \tag{1.9}$$

In addition, under the condition $s_{k-1}^T y_{k-1} > 0$, the two values always satisfy $\alpha_k^{BB_1} \geq \alpha_k^{BB_2}$.

Raydan[2]proved that the BB algorithm with $f(x)$ being strictly convex function is globally convergent, Dai and Liao[3] showed the BB algorithm is of R-linear convergence. For a continuously differentiable objective function, the step sizes $\alpha_k^{BB_1}$ and $\alpha_k^{BB_2}$ of BB algorithm might be negative, and monotonic decreasing property of the objective function is not guaranteed. To overcome the weakness of BB algorithm, Raydan[4] provided with a security strategy on step size and proposed a globally convergent and non-monotonic line search BB algorithm to solve a generic unconstrained optimization problem.

The merits of BB algorithm on computation attracts many researchers' interests and a lot of significant results[5,6,7,8,9]were then achieved. The merits for BB algorithm are based on the fact that construction of step sizes $\alpha_k^{BB_1}$ and $\alpha_k^{BB_2}$ approximately satisfies quasi Newton equations.

However, the determination of $\alpha_k^{BB_1}$ and $\alpha_k^{BB_2}$ only uses the information of the two points $x_k$ and $x_{k-1}$, that is $s_{k-1} = x_k - x_{k-1}$ and $y_{k-1} = g_k - g_{k-1}$. The amount of information is less and couldn't show the advantage of quasi-Newton algorithm.

Inspired by the research work in [1-9], in this paper we use information from the three points $x_k$, $x_{k-1}$ and $x_{k-2}$, or $s_{k-1} = x_k - x_{k-1}$, $s_{k-2} = x_k - x_{k-2}$, $y_{k-1} = g_k - g_{k-1}$ and $y_{k-2} = g_k - g_{k-2}$, and the quasi Newton condition in the sense of least squares to construct new formula for computation of three-point step size gradient.

By using the idea of [4] to secure step size, we design a step size guarantee strategy for the algorithm. Sun Q Y[10]gave a kind of generalized Armijo line search method, i.e. For any real numbers $\mu_1, \mu_2 \in (0,1)$ and $\mu_1 \leq \mu_2$, $\gamma_1, \gamma_2 > 0$, the step size $\lambda_k$ is chosen so that

$$f(x_k + \lambda_k d_k) \leq f(x_k) + \mu_1 \lambda_k \nabla f(x_k)^T d_k, \tag{1.10a}$$

and
$$\lambda_k \geq \gamma_1 \text{ or } \lambda_k \geq \gamma_2 \cdot \lambda_k^* > 0, \tag{1.10b}$$

where $\lambda_k^*$ satisfies:

$$f(x_k + \lambda_k^* d_k) > f(x_k) + \mu_2 \lambda_k^* \nabla f(x_k)^T d_k. \tag{1.10c}$$

Shi Z. J [11] gave a kind of new Armijo line search method with non-exact relaxation, i.e. for $s_k = -(g_k^T d_k)/(L_k \|d_k\|^2)$, take step size $\alpha_k$ to be the maximal in $\{s_k, \omega s_k, \omega^2 s_k, \cdots\}$ satisfying that

$$f(x_k + \alpha d_k) \leq f(x_k) + \sigma\alpha[g_k^T d_k + \frac{1}{2}\alpha\mu L_k \|d_k\|^2]. \tag{1.11}$$

In [11] the proofs of the convergence of algorithm require that the gradient of the objective function $\nabla f(x)$ is Lipschitz continuous. This strong assumption on the objective function restricts application of the algorithms.

Inspired by the research work in [10, 11], in this paper we design relaxed generalized Armijo line search rule and propose three-point step size gradient algorithm with relaxation to solve unconstrained optimization problem, and also we prove that the algorithm is of global convergence under uniform continuity of objective function. In order to test the new algorithms and the analysis on them, numerical experiments are carried out. Computational results confirm the effectiveness of the new algorithm.

In Section 2, we construct the formula of step size for three-point step size gradient algorithm. We present a new method in Section 3. We start the convergence analysis of the new method in Section 4. The convergence properties for generalized convex functions are discussed in Section 5. Finally, a detailed list of the test problems that we have used is given in Section 6.

## 2. Construction of formula of step size for three-point step size gradient algorithm

Before presenting the three-point gradient algorithms, in this section we first introduce how to construct formulas of step size within the algorithms.

For an iterative sequence $\{x_k\}$ produced by an algorithm, objective function $f(x)$ is approximated by the Taylor expansion of second order accuracy at the point $x_k$

$$f(x) \approx f(x_k) + g_k^T(x - x_k) + \frac{1}{2}(x - x_k)^T \nabla^2 f(x_k)(x - x_k),$$

where and what follows $g_k = \nabla f(x_k), B_k = \nabla^2 f(x_k)$. Using Taylor expansion again the gradient of function is approximated by

$$\nabla f(x) \approx g_k + B_k(x - x_k). \tag{2.1}$$

Setting $x = x_{k-1}$ (2.1) yields $g_{k-1} \approx g_k + B_k(x_{k-1} - x_k)$ or equivalently

$$B_k(x_k - x_{k-1}) \approx g_k - g_{k-1}.$$

Letting $s_{k-1} = x_k - x_{k-1}$, $y_{k-1} = g_k - g_{k-1}$, and changing the approximate equal sign into equal sign, we obtain

$$B_k s_{k-1} = y_{k-1} \tag{2.2}$$

or

$$H_k y_{k-1} = s_{k-1}. \tag{2.3}$$

where $B_k = H_k^{-1}$.

Inserting $x = x_{k-2}$ in (2.1), it follows that $g_{k-2} \approx g_k + B_k(x_{k-2} - x_k)$, or

$$B_k(x_k - x_{k-2}) \approx g_k - g_{k-2}.$$

Similarly, taking $s_{k-2} = x_k - x_{k-2}$ and $y_{k-2} = g_k - g_{k-2}$, and changing the signs, we have

$$B_k s_{k-2} = y_{k-2}. \tag{2.4}$$

or

$$H_k y_{k-2} = s_{k-2}. \tag{2.5}$$

By the idea of [1] that $H_k = \alpha_k I$ is supposed to be a scalar matrix, and the purpose to make $H_k$ contain the second order information of objective function, the step size $\alpha_k$ is determined by making it satisfy in the sense of least squares the quasi Newton equation

$$\min_{\alpha>0} \frac{1}{2}\left\|y_{k-1} - \frac{1}{\alpha}s_{k-1}\right\|^2 + \frac{1}{2}\left\|y_{k-2} - \frac{1}{\alpha}s_{k-2}\right\|^2 \tag{2.6}$$

and

$$\min_{\alpha>0} \frac{1}{2}\left\|\alpha y_{k-1} - s_{k-1}\right\|^2 + \frac{1}{2}\left\|\alpha y_{k-2} - s_{k-2}\right\|^2. \tag{2.7}$$

By solving the two quasi Newton equations, a formula to compute step sizes for three-point step size gradient method is obtained

$$\alpha_k^{TBB_1} = \frac{\|s_{k-1}\|^2 + \|s_{k-2}\|^2}{s_{k-1}^T y_{k-1} + s_{k-2}^T y_{k-2}} \tag{2.8}$$

and

$$\alpha_k^{TBB_2} = \frac{s_{k-1}^T y_{k-1} + s_{k-2}^T y_{k-2}}{\|y_{k-1}\|^2 + \|y_{k-2}\|^2}. \tag{2.9}$$

**Remark 2.1.** The formulas (2.8) and (2.9) above are generalization of the ones in the two-point step size gradient algorithm proposed in [1]. In fact, when $x_{k-2} = x_{k-1}$, (2.8) and (2.9) are reduced into the BB formulas (1.8) and (1.9). In addition, formula (2.8) and (2.9) can be extended to the cases with $m$ points.

In order to make full use of information of the three points, combined with classical BB formula $\alpha_k^{BB_1}, \alpha_k^{BB_2}$, let

$$\alpha_k'^{TBB_1} = \lambda \alpha_k^{BB_1} + (1-\lambda)\alpha_k^{BB_2}. \tag{2.10}$$

In order to determine $\lambda$, the quasi-Newton equation (2.4) is used to construct a mathematical model

$$\min_{\lambda \in R} \frac{1}{2}\left\|\frac{1}{\lambda \alpha_k^{BB_1} + (1-\lambda)\alpha_k^{BB_2}} s_{k-2} - y_{k-2}\right\|^2.$$

Thus we have

$$\lambda = \frac{\frac{\|s_{k-2}\|^2}{y_{k-2}^T s_{k-2}} - \alpha_k^{BB_2}}{\alpha_k^{BB_1} - \alpha_k^{BB_2}}, \tag{2.11}$$

where $\alpha_k^{BB_1} = \frac{\|s_{k-1}\|^2}{s_{k-1}^T y_{k-1}}$ and $\alpha_k^{BB_2} = \frac{s_{k-1}^T y_{k-1}}{\|y_{k-1}\|^2}$.

In the same way, let

$$\alpha_k'^{TBB_2} = \lambda \alpha_k^{BB_1} + (1-\lambda)\alpha_k^{BB_2}. \tag{2.12}$$

By using the quasi-Newton equation (2.5), there is

$$\min_{\lambda \in R} \frac{1}{2} \left\| (\lambda \alpha_k^{BB_1} + (1-\lambda) \alpha_k^{BB_2}) y_{k-2} - s_{k-2} \right\|^2.$$

Thus

$$\lambda = \frac{\dfrac{y_{k-2}^T s_{k-2}}{\|y_{k-2}\|^2} - \alpha_k^{BB_2}}{\alpha_k^{BB_1} - \alpha_k^{BB_2}}. \quad (2.13)$$

where $\alpha_k^{BB_1} = \dfrac{\|s_{k-1}\|^2}{s_{k-1}^T y_{k-1}}$ and $\alpha_k^{BB_2} = \dfrac{s_{k-1}^T y_{k-1}}{\|y_{k-1}\|^2}$.

### 3. Algorithm and Properties

With the help of the step size formula constructed in Section 2, in this section we use this formula to present the new algorithm and related properties.

By (2.3), $y_{k-1}^T H_k y_{k-1} = y_{k-1}^T s_{k-1}$. If $H_k = \beta I_n$, then

$$\beta = \frac{s_{k-1}^T y_{k-1}}{\|y_{k-1}\|^2}, \quad (3.1)$$

and $\alpha_{\min} \leq \beta \leq \alpha_{\max}$, where $\alpha_{\min}$ and $\alpha_{\max}$ are the minimum and maximum eigenvalues of $H_k$, respectively.

The new algorithm, named as : Three-pint step size gradient method with relaxed generalized Armijo step size rule (TBB Algorithm), consists of the following steps:

**Step 1** $\forall x_{-1} = x_0 \in R^n$, $L_0 > 0$, $\mu_1, \mu_2 \in (0,1)$ and $\mu_1 \leq \mu_2$, $\gamma_1, \gamma_2 > 0$, $\sigma_1 \in (0,1)$, $\sigma_2 \in (1,2)$, $0 < \alpha_{\min} < \alpha_{\max}$, $0 < \omega < 1$, $\bar{\alpha}_0^{TBB} = 1$, $\varepsilon > 0$. Set $k = 0$.

**Step 2** If $\|g_k\| < \varepsilon$, then stop; Otherwise, go to Step 3.

**Step 3** Let $d_k = -\bar{\alpha}_k^{TBB} g_k$, go to Step 4.

**Step 4** Let $x_{k+1} = x_k + \lambda_k d_k$, where $\lambda_k = \omega^{p_k}$, and $p_k$ is the minimum non-negative integer $p$ satisfying the following formula,

$$f(x_k + \lambda_k d_k) \leq f(x_k) + \mu_1 \lambda_k (\nabla f(x_k)^T d_k + \frac{1}{2\bar{\alpha}_k^{TBB}} L_k \|d_k\|^2), \quad (3.2a)$$

and
$$\lambda_k \geq \gamma_1 \text{ or } \lambda_k \geq \gamma_2 \cdot \lambda_k^* > 0, \quad (3.2b)$$

where $\lambda_k^*$ satisfies:

$$f(x_k + \lambda_k^* d_k) > f(x_k) + \mu_2 \lambda_k^* (\nabla f(x_k)^T d_k + \frac{1}{2\bar{\alpha}_k^{TBB}} L_k \|d_k\|^2) . \tag{3.2c}$$

Where $0 \leq L_k \leq \bar{\alpha}_k^{TBB}(-g_k^T d_k)/\|d_k\|^2$.

**Step 5** Let $k := k+1$, by (2.8) or (2.9) or (2.10) or (2.12) we compute $\alpha_k^{TBB}$, and by (3.1) compute $\beta$.

Denote by

$$\bar{\alpha}_{\min} = \max\{\alpha_{\min}, \sigma_1 |\beta|\}, \quad \bar{\alpha}_{\max} = \min\{\alpha_{\max}, \sigma_2 |\beta|\}.$$

If $\bar{\alpha}_{\min} \leq \alpha_k^{TBB} \leq \bar{\alpha}_{\max}$, choose $\bar{\alpha}_k^{TBB} = \alpha_k^{TBB}$; If $\alpha_k^{TBB} < \bar{\alpha}_{\min}$, choose $\bar{\alpha}_k^{TBB} = \bar{\alpha}_{\min}$;

If $\alpha_k^{TBB} > \bar{\alpha}_{\max}$, take $\bar{\alpha}_k^{TBB} = \bar{\alpha}_{\max}$. Go to Step 2.

**Remark 3.1.** The algorithm above with $\alpha_k^{TBB} = \alpha_k^{TBB_1}$ is denoted by TBB1; that with $\alpha_k^{TBB} = \alpha_k^{TBB_2}$ by TBB2; that with $\alpha_k^{TBB} = \alpha_k^{BB_1}$ is denoted by BB1; that with $\alpha_k^{TBB} = \alpha_k^{BB_2}$ by BB2; that $\alpha_k^{TBB} = \alpha_k'^{TBB_1}$ is denoted by $TBB1'$; that with $\alpha_k^{TBB} = \alpha_k'^{TBB_2}$ by $TBB2'$.

**Remark 3.2.** If $L_k = 0$, the relaxed generalized Armijo step size rule degenerates into the generalized Armijo step size rule. Obviously, if $g_k^T d_k < 0$ and $\lambda_k$ satisfies the generalized Armijo step size rule, $\lambda_k$ must satisfy the relaxed generalized Armijo step size rule.

**Lemma 1** If $x_k$ is not a stationary point of the problem (1.1), then it holds that

$$\|d_k\| \leq \alpha_{\max} \|g_k\|.$$

**Proof.** By the virtue of the definition of $d_k$ it is easy to prove, here omitted.

**Lemma 2** Assume that $x_k$ is not a stationary point of the optimization problem (1.1), then it holds that $g_k^T d_k \leq -\alpha_{\min} \|g_k\|^2$.

**Proof.** By the definition of $d_k$ this is easily proved.

## 4. Global convergence of algorithm

Throughout this paper, let $\{x_k\}$ be a sequence of points produced by TBB (one of TBB1, TBB2, BB1, BB2, $TBB1'$, $TBB2'$). If there exists an index $k$ such that $\nabla f(x_k) = 0$, then $x_k$ is a stationary point for the optimization problem (1.1). In what follows, $\{x_k\}$ produced by TBB is assumed to be an infinite sequence. In this section the global convergence of $\{x_k\}$ is proved.

**Theorem 1** suppose that $f(x) \in C^1$. Then:

（i） either $f(x_k) \to -\infty (k \to \infty)$, or $\liminf\limits_{k \to \infty} \|\nabla f(x_k)\| = 0$;

（ii）either $f(x_k) \to -\infty (k \to \infty)$, or $\lim\limits_{k \to \infty} \|\nabla f(x_k)\| = 0$, if $\nabla f(x)$ is uniformly continuous on $R^n$.

**Proof.** By Lemma2, for all $k$, we have that $\nabla f(x_k)^T d_k < 0$. It follows from（3.1a）that

$$f(x_k + \lambda_k d_k) \le f(x_k) + \mu_1 \lambda_k (\nabla f(x_k)^T d_k + \frac{1}{2\overline{\alpha}_k^{TBB}} L_k \|d_k\|^2)$$

$$\le f(x_k) + \frac{1}{2} \mu_1 \lambda_k \nabla f(x_k)^T d_k$$

$$\le f(x_k),$$

which implies that $\{f(x_k)\}$ is a monotonically decreasing sequence. If $f(x_k) \to -\infty (k \to \infty)$, then we complete the proof. Therefore, in the following discussion, we assume that $\{f(x_k)\}$ is a bounded set.

Suppose （i） is not true. Then, there exists $\varepsilon > 0$ such that

$$\|\nabla f(x_k)\| \ge \varepsilon, \ \forall k. \tag{4.1}$$

It follows from Lemma 2 and（3.2a）that

$$f(x_{k+1}) - f(x_k) \le \frac{1}{2} \mu_1 \lambda_k \nabla f(x_k)^T d_k \le -\frac{1}{2} \alpha_{\min} \cdot \lambda_k \mu_1 \|\nabla f(x_k)\|^2. \tag{4.2}$$

The above inequality（4.1），（4.2）and the boundedness of $\{f(x_k)\}$ imply that

$$\sum_{k=1}^{\infty} \lambda_k \|\nabla f(x_k)\| < +\infty. \tag{4.3}$$

It follows from Lemma 1 that

$$\|x_{k+1} - x_k\| = \lambda_k \|d_k\| \le \lambda_k \alpha_{\max} \|\nabla f(x_k)\|.$$

The above inequalities and (4.3) yield $\sum_{k=1}^{\infty}\|x_{k+1} - x_k\| < +\infty$, which yields that $\{x_k\}$ is convergent, say to a pint $x_*$. From (4.1) and (4.3), we have

$$\lim_{k \to \infty} \lambda_k = 0. \qquad (4.4)$$

It follows from Lemma1, the convergence of $\{x_k\}$ and $f(x) \in C^1$ that $\{d_k\}$ is bounded. Without loss of generality, we may assume that there exists an index set $k \subset \{1,2,...\}$ such that $\lim_{k \to \infty, k \in K} d_k = d_*$. It follows from (4.4) and (3.2b) that, when $k (k \in K)$ is large enough, we have $\lambda_k < \gamma_1$, and hence it follows from (3.2b) that, $\lambda_k \geq \gamma_2 . \lambda_k^*$, where $\lambda_k^*$ satisfies (3.2c), i.e.

$$(f(x_k + \lambda_k^* d_k) - f(x_k))/\lambda_k^* > \mu_2 \nabla f(x_k)^T d_k.$$

Taking the limit for $k \in K$, we have

$$\nabla f(x_*)^T d_* \geq \mu_2 \nabla f(x_*)^T d_*.$$

By using the above inequality and $\mu_2 \in (0,1)$, we obtain that

$$\nabla f(x_*)^T d_* = 0.$$

The above equality and Lemma 2 yield $\|\nabla f(x_*)\| = 0$, which contradicts (4.1). This completes the proof of (i).

Suppose that $\lim_{k \to \infty} \|\nabla f(x_k)\| = 0$ is not true. Then, there exist an infinite index set $K_1 \subset \{1,2,...\}$ and a positive scalar $\varepsilon > 0$ such that, for all $k \in K_1$,

$$\|\nabla f(x_k)\| > \varepsilon. \qquad (4.5)$$

It follows from Lemma 2 and (3.2a) that

$$f(x_k) - f(x_{k+1}) \geq -\frac{1}{2}\mu_1 \lambda_k \nabla f(x_k)^T d_k \geq \frac{1}{2}\alpha_{\min} \lambda_k \mu_1 \|\nabla f(x_k)\|^2. \qquad (4.6)$$

By using (4.5) and (4.6), we obtain that

$$\lambda_k \leq \frac{2}{\mu_1 \varepsilon^2} \alpha_{\min}^{-1} (f(x_k) - f(x_{k+1})), \forall k \in K_1.$$

The boundedness of $\{f(x_k)\}$ and the monotonically decreasing property imply that $\{f(x_k)\}$ is convergent. Thus,

$$\limsup_{k\to\infty, k\in K_1} \lambda_k \le \limsup_{k\to\infty, k\in K_1} \frac{2}{\mu_1 \varepsilon^2} \alpha_{\min}^{-1} (f(x_k) - f(x_{k+1})),$$

Which yields that

$$\lim_{k\to\infty, k\in K_1} \lambda_k = 0. \tag{4.7}$$

It follows from (4.5) and (4.6) that

$$\lambda_k \|\nabla f(x_k)\| \le \frac{2}{\mu_1 \varepsilon} \alpha_{\min}^{-1} (f(x_k) - f(x_{k+1})),$$

and

$$\limsup_{k\to\infty, k\in K_1} \lambda_k \|\nabla f(x_k)\| \le \limsup_{k\to\infty, k\in K_1} \frac{2}{\mu_1 \varepsilon} \alpha_{\min}^{-1} (f(x_k) - f(x_{k+1})).$$

Hence,

$$\limsup_{k\to\infty, k\in K_1} \lambda_k \|\nabla f(x_k)\| = 0. \tag{4.8}$$

It follows from Lemma 1 and (4.8) that

$$\limsup_{k\to\infty, k\in K_1} \lambda_k \|d_k\| \le \limsup_{k\to\infty, k\in K_1} \alpha_{\max} \lambda_k \|\nabla f(x_k)\| = 0.$$

i.e.

$$\limsup_{k\to\infty, k\in K_1} \lambda_k \|d_k\| = 0. \tag{4.9}$$

It follows from (4.7) that, when $k(k \in K_1)$ is large enough, we have $\lambda_k < \gamma_1$, and hence it follows from (3.2b) that, $\lambda_k \ge \gamma_2 . \lambda_k^*$, where $\lambda_k^*$ satisfies (3.2c). Now set $x_{k+1}^* = x_k + \lambda_k^* d_k$. It follows from (4.7) and (4.9) and $\lambda_k \ge \gamma_2 . \lambda_k^*$ ($k \in K_1$ is large enough) that $\lim_{k\to\infty, k\in K_1} \lambda_k^* = 0$ and $\lim_{k\to\infty, k\in K_1} \lambda_k^* \|d_k\| = 0$. Hence, $\lim_{k\to\infty, k\in K_1} \|x_{k+1}^* - x_k\| = 0$.

Let $\rho_k^* = \dfrac{f(x_{k+1}^*) - f(x_k)}{\lambda_k^* \nabla f(x_k)^T d_k}, k \in K_1$, it follows from (3.2c) that

$$\rho_k^* < \mu_2 < 1, k \in K_1. \tag{4.10}$$

It follows from Lemma 1, 2 and (4.5) that

$$\limsup_{k\to\infty, k\in K_1} |\rho_k^* - 1| = \limsup_{k\to\infty, k\in K_1} \left| \frac{\nabla f(\xi_k^*)^T (\lambda_k^* d_k)}{\lambda_k^* \nabla f(x_k)^T d_k} - 1 \right|$$

$$= \limsup_{k\to\infty, k\in K_1} \left| \frac{(\nabla f(\xi_k^*) - \nabla f(x_k))^T d_k}{\nabla f(x_k)^T d_k} \right|$$

$$\leq \limsup_{k\to\infty, k\in K_1} \frac{\|\nabla f(\xi_k^*) - \nabla f(x_k)\| \cdot \|d_k\|}{|\nabla f(x_k)^T d_k|}$$

$$\leq \limsup_{k\to\infty, k\in K_1} \frac{\|\nabla f(\xi_k^*) - \nabla f(x_k)\| \cdot \alpha_{max} \cdot \|\nabla f(x_k)\|}{\alpha_{min} \cdot \|\nabla f(x_k)\|^2}$$

$$\leq \limsup_{k\to\infty, k\in K_1} \frac{\|\nabla f(\xi_k^*) - \nabla f(x_k)\| \cdot \alpha_{max} \cdot \|\nabla f(x_k)\|}{\alpha_{min} \cdot \varepsilon^2} = 0. \qquad (4.11)$$

Where $\xi_k^* = x_k + \vartheta_k(x_{k+1}^* - x_k), 0 < \vartheta_k < 1, k \in K_1$.

Hence (4.11) establishes that $\rho_k^* = 1$. This is the desired contradiction because (4.10) guarantees that $\rho_k^* < \mu_2 < 1$. This yields (ii).

**Remark 4.1** Following the proof of Theorem 3.1 in [13], it is easy to prove that the TBB algorithm is linearly convergent. Following the proof of Theorem 2 in [12], it is easy to prove that the algorithm TBB is super-linearly convergent.

## 5. Convergence Properties for Generalized Convex Functions

In this section we discuss the convergence properties of (TBB) for generalized convex functions. In the following discussion, we make the following assumption:

**Assumption 5.1**: $\lambda_0 = \sup\{\lambda_k, k = 1, 2, ...\} < +\infty$.

As shown in the following, Lemma 3 play an important role in our analysis. Thus we have the following results.

**Lemma 3** Suppose that Assumption 5.1 holds and $f(x) \in C^1$. If $f(x)$ is a quasi-convex function and the solution set of problem (1.1) is nonempty, then $\{x_k\}$ is a bounded sequence, and each accumulation point $x_*$ of which is a stationary point of problem (1.1) and $\lim_{k\to\infty} x_k = x_*$.

**Proof.** Note that for all $x \in R^n$ and all $k$,

$$\|x_{k+1} - x\|^2 = \|x_k - x\|^2 + 2(x_{k+1} - x_k, x_k - x) + \|x_{k+1} - x_k\|^2$$

$$= \|x_k - x\|^2 + 2\lambda_k(d_k, x_k - x) + \lambda_k^2 \|d_k\|^2$$

$$= \|x_k - x\|^2 + 2\lambda_k \bar{\alpha}_k^{TBB}(-\nabla f(x_k), x_k - x) + \lambda_k^2 \|d_k\|^2. \qquad (5.1)$$

It follows from Lemma 1, Lemma 2 and (3.2a) that

$$\|d_k\|^2 \leq \alpha_{max}^2 \|\nabla f(x_k)\|^2; \qquad (5.2)$$

$$\|\nabla f(x_k)\|^2 \leq \alpha_{min}^{-1} \cdot (-\nabla f(x_k)^T d_k); \qquad (5.3)$$

$$-\lambda_k \nabla f(x_k)^T d_k \leq 2\mu_1^{-1}(f(x_k) - f(x_{k+1})) .\tag{5.4}$$

By using (5.1), (5.2), (5.3) and (5.4), we obtain that

$$\|x_{k+1} - x\| \leq \|x_k - x\|^2 + 2\lambda_0 (\nabla f(x_k), x - x_k) + \lambda_0 \lambda_k \alpha_{max}^2 \cdot \alpha_{min}^{-1} \cdot (-\nabla f(x_k)^T d_k)$$

$$\leq \|x_k - x\|^2 + 2\lambda_0 (\nabla f(x_k), x - x_k) + \lambda_0 \alpha_{max}^2 \cdot \alpha_{min}^{-1} \mu_1^{-1}(f(x_k) - f(x_{k+1}))$$

$$= \|x_k - x\|^2 + 2\lambda_0 (\nabla f(x_k), x - x_k) + m(f(x_k) - f(x_{k+1})) ,\tag{5.5}$$

where $m = \lambda_0 \mu_1^{-1} \alpha_{max}^2 \cdot \alpha_{min}^{-1}$.

Because the solution set of problem (1.1) is nonempty, we can choose $y \in R^n$ satisfying $f(y) \leq f(x_k)$. Since $f(x)$ is a quasi-convex function, we have

$$(\nabla f(x_k), y - x_k) \leq 0 .\tag{5.6}$$

It follows from (5.5), (5.6) that

$$\|x_{k+1} - y\|^2 + mf(x_{k+1}) \leq \|x_k - y\|^2 + mf(x_k) ,$$

Which implies the sequence $\{\|x_k - y\|^2 + mf(x_k)\}$ is descent. Since we have assumed that the solution set of problem (1.1) is nonempty, and so $\inf\{f(x_k) : k = 1,2,...\} > -\infty$, both sequence $\{f(x_k)\}$ and $\{\|x_k - y\|^2 + mf(x_k)\}$ are bounded from below and converge. Therefore, the sequence $\{\|x_k - y\|^2\}$ converges and $\{x_k\}$ is bounded. This implies that $\{x_k\}$ has an accumulation point $x_*$ and that there exists an index set $K_1 \subset \{1,2,...\}$ such that

$$\lim_{k \to \infty, k \in K_1} x_k = x_*, \quad \text{and} \quad \lim_{k \to \infty, k \in K_1} f(x_k) = f(x_*).$$

It follows from the above equation and the fact $\{f(x_k)\}$ is a monotonically decreasing sequence implies $\lim_{k \to \infty, k \in K_1} f(x_{k-1}) = f(x_*)$. Therefore, we have

$$\lim_{k \to \infty} \{ \|x_k - x_*\|^2 + mf(x_k) \}$$

$$= \lim_{k \to \infty, k \in K_1} \{ \|x_k - x_*\|^2 + mf(x_k) \}$$

$$= mf(x_*),$$

Which implies $\lim_{k \to \infty} x_k = x_*$. From Theorem 1 the limit $x_*$ is a stationary point of problem (1.1)

**Theorem 2** Suppose that Assumption 5.1 holds and $f(x) \in C^1$. If $f(x)$ is a pseudo-convex function, then:

(i) $\{x_k\}$ is abounded sequence if and only if the solution set of problem (1.1) is nonempty;

(ii) $\lim_{k \to \infty} f(x_k) = \inf\{f(x) : x \in R^n\}$;

(iii) If the solution set of problem (1.1) is nonempty, then any accumulation point $x_*$ of $\{x_k\}$ is an optimal solution of problem (1.1) and $\lim_{k \to \infty} x_k = x_*$.

**Proof.** Since $f(x)$ is pseudo-convex, it is quasi-convex and a stationary point of problem (1.1) is also an optimal solution of problem (1.1).

First, we will show part (i). If $\{x_k\}$ is a bounded sequence, then it follows from Theorem 1 that there exists an index set $K_2 \subset \{1,2,...\}$ and a point $x_* \in R^n$ such that $\lim_{k \to \infty, k \in K_2} x_k = x_*$, and $x_*$ is a stationary point of problem (1.1), and is also an optimal solution of problem (1.1). Conversely, if the solution set of problem (1.1) is nonempty, then it follows from Lemma 3 that $\{x_k\}$ is a bounded sequence.

Next, we will prove (ii). We prove this conclusion by the following three cases (a), (b), (c).

(a) $\lim_{k \to \infty} f(x_k) = \inf\{f(x_k) : k = 1,2,...\} = -\infty$; It follows from $\{f(x_k)\}$ is a descent sequence, and $\lim_{k \to \infty} f(x_k) = \inf\{f(x_k) : k = 1,2,...\} \geq \inf\{f(x) : x \in R^n\}$.

(b) $\{x_k\}$ is bounded: It follows from (i) of this theorem that the solution set of problem (1.1) is nonempty, and there exists an index set $K_3 \subset \{1,2,...\}$ and a point $x_* \in R^n$ such that $\lim_{k \to \infty, k \in K_3} x_k = x_*$, it follows from Theorem 1 that $x_*$ is a stationary point of problem (1.1), and is also an optimal solution of problem (1.1).

(c) $\inf\{f(x_k) : k = 1,2,...\} > -\infty$, and $\{x_k\}$ is unbounded: suppose that there exists $\bar{x} \in R^n$, $\varepsilon > 0$, and $k_1$ such that for all $k \geq k_1$, $f(x_k) > f(\bar{x}) + \varepsilon$. Since $f(x)$ is a pseudo-convex function, we have $(\nabla f(x_k), \bar{x} - x_k) \leq 0$, for all $k \geq k_1$. Setting $x = \bar{x}$ in (5.5) that

$$\|x_{k+1} - \bar{x}\|^2 + mf(x_{k+1}) \leq \|x_k - \bar{x}\|^2 + mf(x_k),$$

Which implies the sequence $\{\|x_k - \bar{x}\|^2 + mf(x_k)\}$ is descent. Since we assumed that $\inf\{f(x_k) : k = 1,2,...\} > -\infty$; both sequence $\{f(x_k)\}$ and $\{\|x_k - \bar{x}\|^2 + mf(x_k)\}$ are bounded from below and converge. Therefore, the sequence $\{\|x_k - \bar{x}\|^2\}$ converges and $\{x_k\}$ is bounded, which contradicts our assumption.

**Corollary 1** Suppose that Assumption 5.1 holds and $f(x) \in C^1$. If $f(x)$ is a convex function, then:

（i） $\{x_k\}$ is abounded sequence if and only if the solution set of problem (1.1) is nonempty;

（ii） $\lim_{k \to \infty} f(x_k) = \inf\{f(x) : x \in R^n\}$;

（iii） If the solution set of problem (1.1) is nonempty, then any accumulation point $x_*$ of $\{x_k\}$ is an optimal solution of problem （1.1） and $\lim_{k \to \infty} x_k = x_*$.

**Proof.**  Since $f(x)$ is convex , it is pseudo-convex. It immediately follows from Theorem 2.

**Corollary 2** Suppose that Assumption 5.1 holds and $f(x) \in C^1$.  If $f(x)$ is a quasi-convex function, then either the solution set of problem (1.1)is empty or any accumulation point $x_*$ of $\{x_k\}$ is an optimal solution of problem （1.1） and $\lim_{k \to \infty} x_k = x_*$.

**Proof.**  It immediately follows from Lemma 3.

## 6. Linear convergence

**Assumption 6.1** The level set $L(x_0) = \{x \in R^n \mid f(x) \leq f(x_0)\}$ of $f \in C^1$ is bounded.

**Assumption 6.2** Let $f(x)$ be quadratic continuous differentiable in the neighborhood of minimum point $x^*$, and there exist $\varepsilon > 0$ and $M > m > 0$, so that when $\|x - x^*\| < \varepsilon$, there is

$$m\|y\|^2 \leq y^T \nabla^2 f(x) y \leq M\|y\|^2 \quad \text{for} \quad \forall y \in R^n.$$

**Lemma 4** Suppose that Assumption 6.2 holds. When $\|x - x^*\| < \varepsilon$, we have

$$\frac{1}{2}m\|x - x^*\|^2 \leq f(x) - f(x^*) \leq \frac{1}{2}M\|x - x^*\|^2, \tag{6.1}$$

$$\|g(x)\| \geq m\|x - x^*\|. \tag{6.2}$$

**Proof** By $\nabla f(x^*) = 0$ and Taylor expansion, we obtain that for $\forall x \in N(x^*, \sigma)$, there exists $\xi \in (x, x^*)$ satisfying

$$f(x) - f(x^*) = \frac{1}{2}(x - x^*)^T \nabla^2 f(\xi)(x - x^*) \ .$$

By Assumption 6.2, (6.1) is proved.

Combining

$$\nabla f(x) = \nabla f(x) - \nabla f(x^*) = \int_0^1 \nabla^2 f(x^* + \tau(x - x^*))(x - x^*) d\tau$$

with Assumption 6.2, we have

$$\|x - x^*\| \cdot \|\nabla f(x)\| \geq (x - x^*)^T \nabla f(x)$$

$$= \int_0^1 (x - x^*)^T \nabla^2 f(x^* + \tau(x - x^*))(x - x^*) d\tau$$

$$\geq m \|x - x^*\|^2 .$$

Thus $\|g(x)\| \geq m \|x - x^*\|$.

**Lemma 5** Suppose that Assumption 6.1 and 6.2 hold, and the sequence $\{x_k\}$ generated by the algorithm converges to $x^*$, then $\inf \lambda_k = \eta_0 > 0$.

**Proof** If the algorithm terminates in limited steps, the conclusion is true.

Supposing that the sequence $\{x_k\}$ is infinite. As $0 < \lambda_k \leq 1$, $\inf \lambda_k = \eta_0 \geq 0$. Next, we show that $\eta_0 = 0$ is impossible.

We prove by contradiction. Suppose that $\eta_0 = 0$, then there exists a subsequence $\{\lambda_k\}_{k \in K_1}$ such that

$$\lambda_k \to 0, k \in K_1. \qquad (6.3)$$

As

$$\|d_k\| = \overline{\alpha}_k^{TBB} \|g_k\| 及 \ x_k \to x^*, g_k \to g^*, k \to \infty,$$

and by Assumption 6.1 and $\overline{\alpha}_{\min} \leq \overline{\alpha}_k^{TBB} \leq \overline{\alpha}_{\max}$, $\|d_k\|$ is bounded on $L(x_0)$. So there exists a constant $c$ such that $\|d_k\| \leq c$. By (6.3), when $k \in K_1$ is sufficiently large, we have $\lambda_k \geq \gamma_2 . \lambda_k^*$. Therefore,

$$\lim_{k \to \infty, k \in K_1} \lambda_k^* = 0 \ \text{and} \ \lim_{k \to \infty, k \in K_1} \lambda_k^* \|d_k\| = 0,$$

that is
$$\lim_{k\to\infty, k\in K_1} \|x^*_{k+1} - x_k\| = 0.$$

There is a $\bar{k}$ satisfied that $\|\lambda_k^* d_k\| < \frac{\varepsilon}{2}$ and $\|x_k - x^*\| < \frac{\varepsilon}{2}$ hold for $k > \bar{k}$ ($k \in K_1$). Thus, when $k > \bar{k}$ ($k \in K_1$), we have

$$\|x_k + \lambda_k^* d_k - x^*\| \leq \|x_k - x^*\| + \|\lambda_k^* d_k\| < \varepsilon. \tag{6.4}$$

By $\lambda_k^*$, there is

$$f(x_k + \lambda_k^* d_k) > f(x_k) + \lambda_k^* \gamma_0 [g_k^T d_k + \frac{1}{2\bar{\alpha}_k^{TBB}} L_k \|d_k\|^2]$$

$$> f(x_k) + \lambda_k^* \gamma_0 g_k^T d_k,$$

that is

$$f(x_k) - f(x_k + \lambda_k^* d_k) < -\lambda_k^* \gamma_0 g_k^T d_k.$$

From the mean value theorem, there is $\theta_k \in [0,1]$ such that

$$f(x_k) - f(x_k + \lambda_k^* d_k) = -\lambda_k^* g(x_k + \theta_k \lambda_k^* d_k)^T d_k.$$

According to the above two formulas, there is

$$[g(x_k + \lambda_k^* \theta_k d_k) - g(x_k)]^T d_k > -(1-\gamma_0) g_k^T d_k. \tag{6.5}$$

Using the mean value theorem of vector function, it obtains

$$(g(x_k + \lambda_k^* \theta_k d_k) - g(x_k))^T d_k = \lambda_k^* \theta_k \int_0^1 d_k^T \nabla^2 f(x_k + t(\lambda_k^* \theta_k d_k)) d_k dt,$$

i.e.,

$$g(x_k + \lambda_k^* \theta_k d_k) - g(x_k) = \lambda_k^* \theta_k \int_0^1 \nabla^2 f(x_k + t(\lambda_k^* \theta_k d_k)) d_k dt.$$

By (6.4), we have

$$\|x_k + t\lambda_k^* \theta_k d_k - x^*\| < \varepsilon, (t \in [0,1], 0 < \lambda_k \leq 1, \theta_k \in (0,1)).$$

By the above formula and Assumption 6.2, there is

$$\lambda_k^* \theta_k \int_0^1 d_k^T \nabla^2 f(x_k + t(\lambda_k^* \theta_k d_k)) d_k dt \leq \lambda_k^* \int_0^1 M \|d_k\|^2 dt = \lambda_k^* M \|d_k\|^2.$$

By (6.5), it follows

$$\lambda_k^* M \|d_k\|^2 > -(1-\gamma_0) g_k^T d_k,$$

that is

$$\lambda_k^* > \frac{-(1-\gamma_0) g_k^T d_k}{M \|d_k\|^2} = \frac{(1-\sigma) d_k^T \bar{\alpha}_k^{TBB} d_k}{M \|d_k\|^2} \geq \frac{(1-\sigma)\alpha_{\min}}{M} > 0, k \in K_1.$$

This contradicts $\{\lambda_k\} \to 0 (k \in K_1)$. Thus, $\eta_0 > 0$.

**Theorem 3** Suppose that Assumption 6.1 and 6.2 hold, the sequence $\{x_k\}$ generated by the algorithm converges to $x^*$, and

$$0 < \alpha_{\min} \leq \bar{\alpha}_k^{TBB} < \frac{M}{2\eta_0 m^2 \mu_1},$$

then $\{x_k\}$ is R-linear convergence to $x^*$.

**Proof** Suppose that $\|x_k - x^*\| < \varepsilon$ for $\forall k \leq k_1$. By the definition of $\lambda_k$, we have

$$f(x_k) - f(x_k + \lambda_k d_k) \geq \mu_1 \lambda_k (\nabla f(x_k)^T d_k + \frac{1}{2\bar{\alpha}_k^{TBB}} L_k \|d_k\|^2)$$

$$\geq -\lambda_k \mu_1 g_k^T d_k$$

$$\geq \eta_0 \mu_1 \alpha_{\min} \|g_k\|^2$$

$$\geq \eta_0 \mu_1 \alpha_{\min} m^2 \|x_k - x^*\|^2$$

$$\geq \frac{2\eta_0 \mu_1 \alpha_{\min} m^2}{M} [f(x_k) - f(x^*)]$$

$$= \delta^2 [f(x_k) - f(x^*)],$$

where $\delta = \sqrt{\frac{2\eta_0 \mu_1 \alpha_{\min} m^2}{M}}$. Hence,

$$f(x_{k+1}) - f(x^*) \leq [1 - \delta^2][f(x_k) - f(x^*)]. \tag{6.6}$$

Set $\theta = [1 - \delta^2], \theta \in (0,1)$, then

$$f(x_k) - f(x^*) \leq \theta^2 [f(x_{k-1}) - f(x^*)] \leq \cdots \leq \theta^{2k} [f(x_0) - f(x^*)].$$

By (6.1), it follows

$$\|x_k - x^*\|^2 \leq \frac{2}{m}[f(x_k) - f(x^*)] \leq \frac{2}{m}[f(x_0) - f(x^*)]\theta^{2k},$$

that is

$$\|x_k - x^*\| \leq \rho \theta^k,$$

where $\theta < 1, \rho = \sqrt{\frac{2}{m}}[f(x_0) - f(x^*)]^{\frac{1}{2}}$. Therefore

$$R_1\{x_k\} = \lim_{k\to\infty}(\|x_k - x^*\|)^{\frac{1}{k}} = \theta < 1.$$

Thus $\{x_k\}$ is R-linear convergence to $x^*$.

## 7. Superlinear convergence

**Lemma 6**[14] Suppose that $F: R^n \to R^m$ is continuously differentiable on open convex set $D$, then for any $x, u, v \in D$, there is

$$\|F(u) - F(v) - F'(x)(u-v)\| \le \left[\sup_{0 \le t \le 1}\|F'(v + t(u-v)) - F'(x)\|\right]\|u-v\|.$$

Furthermore, if $F'$ is Lipschitz continuous on $D$, there are

$$\|F(u) - F(v) - F'(x)(u-v)\| \le r\frac{\|u-x\| + \|v-x\|}{2}\|u-v\|,$$

and

$$\|F(u) - F(v) - F'(x)(u-v)\| \le r\sigma(u,v)\|u-v\|,$$

where $r$ is the Lipschitz constant, and $\sigma(u,v) = \max\{\|u-x\|, \|v-x\|\}$.

**Lemma 7**[14] Suppose that $F: R^n \to R^m$ is continuously differentiable on open convex set $D$, and $F'$ is Lipschitz continuous on $D$. If $[F'(x)]^{-1}$ exists, there are $\varepsilon > 0$ and $\beta > \alpha > 0$ such that for any $\forall u, v \in D$, when $\max\{\|u-x\|, \|v-x\|\} \le \varepsilon$, we have $\alpha\|u-v\| \le \|F(u) - F(v)\| \le \beta\|u-v\|$.

**Theorem 4** Suppose that (a) the sequence $\{x_k\}$ generated by the algorithm converges to local minimum $x^*$, and $x_k \ne x^*$;

(b) the level set $L(x_0) = \{x \in R^n \mid f(x) \le f(x_0)\}$ is bounded;

(c) $f(x)$ is second order continuously differentiable in a neighborhood of $x^*$, and there exist $\varepsilon > 0$ and $M > m > 0$, such that $\|y - x^*\| < \varepsilon$ and $m\|y\|^2 \le y^T G(x) y \le M\|y\|^2$ for $\forall y \in R^n$;

(d) $G(x)$ is Lipschitz continuous in a neighborhood of $x^*$, where $r$ is the Lipschitz

constant.

Then, $\{x_k\}$ superlinearly converges to $x^*$ if and only if

$$\lim_{k\to\infty}\left[\left\|(\frac{1}{\lambda_k\bar{\alpha}_k^{TBB}}I_n - G(x^*))s_k\right\|\Big/\|s_k\|\right] = 0, \tag{7.1}$$

where $s_k = x_{k+1} - x_k$.

**Proof** Firstly, we prove the necessity. By (b), (c) and (d), it follows that there exists an open convex set $D = N(x^*,\varepsilon) \cap L(x_0)$ such that $x_k \in D$ for sufficiently large $k$. By $d_k = -\bar{\alpha}_k^{TBB} g_k$, we have

$$(\frac{1}{\lambda_k\bar{\alpha}_k^{TBB}}I_n - G(x^*))s_k = -g(x_k) - G(x^*)s_k$$

$$= [g(x_{k+1}) - g(x_k) - G(x^*)s_k] - g(x_{k+1}). \tag{7.2}$$

From the mean value theorem, it obtains

$$\frac{\|g(x_{k+1})\|}{\|s_k\|} \leq \frac{\left\|(\frac{1}{\lambda_k\bar{\alpha}_k^{TBB}}I_n - G(x^*))s_k\right\|}{\|s_k\|} + \frac{\|g(x_{k+1}) - g(x_k) - G(x^*)s_k\|}{\|s_k\|}$$

$$= \frac{\left\|(\frac{1}{\lambda_k\bar{\alpha}_k^{TBB}}I_n - G(x^*))s_k\right\|}{\|s_k\|} + \frac{\left\|\int_0^1 [G(x_k + ts_k) - G(x^*)]s_k dt\right\|}{\|s_k\|},$$

where $t \in [0,1]$. By $\lim_{k\to\infty} x_k = x^*$ and the continuity of $G(x)$, we know

$$\lim_{k\to\infty} \frac{\left\|\int_0^1 [G(x_k + ts_k) - G(x^*)]s_k dt\right\|}{\|s_k\|} = 0. \tag{7.3}$$

By (7.3), it follows

$$\lim_{k\to\infty} \frac{\|g(x_{k+1})\|}{\|s_k\|} = 0. \tag{7.4}$$

By (7.4) and $\lim_{k\to\infty}\|s_k\| = 0$, $g(x^*) = \lim_{k\to\infty} g(x_k) = 0$. From (c) it obtains that $G(x^*)$ is nonsingular. By Lemma 7, there exist $\alpha > 0$ and $k_0 > 0$ such that for $\forall k \geq k_0$,

$$\|g(x_{k+1})\| = \|g(x_{k+1}) - g(x^*)\| \geq \alpha\|x_{k+1} - x^*\|. \tag{7.5}$$

By (7.5), we have

$$\frac{\|g(x_{k+1})\|}{\|x_{k+1} - x_k\|} \geq \frac{\alpha\|x_{k+1} - x^*\|}{\|x_{k+1} - x^*\| + \|x_k - x^*\|} \geq \alpha \frac{r_k}{1 + r_k},$$

where $r_k = \frac{\|x_{k+1} - x^*\|}{\|x_k - x^*\|}$. By (7.1) and (7.3), it follows $\lim_{k \to \infty} \frac{r_k}{1 + r_k} = 0$. Accordingly, it deduces $\lim_{k \to \infty} r_k = 0$, and $\{x_k\}$ superlinearly converges to $x^*$.

Secondly, we prove the sufficiency. By Lemma 7, there exist $\beta > 0$ and $k_0 > 0$ such that for $\forall k \geq k_0$,

$$\|g(x_{k+1})\| \leq \beta \|x_{k+1} - x^*\|. \tag{7.6}$$

Because the sequence $\{x_k\}$ superlinearly converges to $x^*$, by (7.6), we know

$$0 = \lim_{k \to \infty} \frac{\|x_{k+1} - x^*\|}{\|x_k - x^*\|} \geq \lim_{k \to \infty} \frac{\|g(x_{k+1})\|}{\beta \|x_k - x^*\|} = \lim_{k \to \infty} \frac{1}{\beta} \frac{\|g(x_{k+1})\|}{\beta \|x_{k+1} - x_k\|} \cdot \frac{\|x_{k+1} - x_k\|}{\|x_k - x^*\|}. \tag{7.7}$$

By (7.7) and $\frac{\|x_{k+1} - x^*\|}{\|x_k - x^*\|} \geq \left|\frac{\|x_{k+1} - x_k\|}{\|x_k - x^*\|} - \frac{\|x_k - x^*\|}{\|x_k - x^*\|}\right|$, it follows $\lim_{k \to \infty} \frac{\|x_{k+1} - x_k\|}{\|x_k - x^*\|} = 1$, and then

$$\lim_{k \to \infty} \frac{\|g(x_{k+1})\|}{\|x_{k+1} - x_k\|} = 0. \tag{7.8}$$

By (7.2) and (7.3), we have

$$\lim_{k \to \infty} \left[\left\|\left(\frac{1}{\lambda_k \bar{\alpha}_k^{TBB}} I_n - G(x^*)\right)s_k\right\| \bigg/ \|s_k\|\right] = 0.$$

## 8. Numerical Experiments

In this section, we implement algorithm TBB and BB. The programs are all written in Matlab and run on personal computers. We carry out 91 test problems from Andrei test function collection [15]. The dimension of each test problem is 10000.

The parameters are taken as follows. $\mu_1 = \mu_2 = 0.32$, $\omega = 0.76$, $\alpha_{\min} = 0.006$, $\alpha_{\max} = 100$, $\sigma_1 = 0.52$, $\sigma_2 = 1.2$, $L_k = \bar{\alpha}_k^{TBB}(-g_k^T d_k)/\|d_k\|^2$. The algorithm will stop when any of the following termination conditions is met, (1) $\|g(x_k)\| < 10^{-6}\|g(x_0)\|$; (2) The iteration times is more than 10000; (3) The CPU time exceeds 600 seconds.

In this paper, the numerical performance comparison method proposed by Dolan et al. [16] is used to give the comparison results of the numerical performance of different algorithms. The performance diagram in [16] is illustrated by taking the iteration times as an example. The function $P(\tau)$ represents the ratio of the number of problems that the test algorithms can solve to the total number of problems within $\tau$ times of the minimum iteration set $I^*$. It reflects the numerical performance of the iteration times of the test algorithm. For each test algorithm, the curve of $P(\tau)$ is given. The highest curve indicates that the test algorithm can solve the most problems within $\tau$ times of the minimum iteration set, which means that the test algorithm has better numerical performance in iteration times. The numerical experiment is divided into two parts.

In the first part, the numerical performance comparison method [16] is used to compare BB1, BB2, TBB1, TBB2, $TBB1'$ and $TBB2'$. Under the termination conditions, $TBB1'$ can successfully solve 60 problems, while other algorithms can successfully solve 54 of them. The following analysis only considers the 60 problems that can be successfully solved by $TBB1'$. Figures 1 to 4 show the performance evaluation results of the 6 algorithms in terms of iteration times, function value calculation times, gradient calculation times and CPU time. It is easy to see that the numerical effects of $TBB1'$ and TBB1 are better than those of BB1 and BB2.

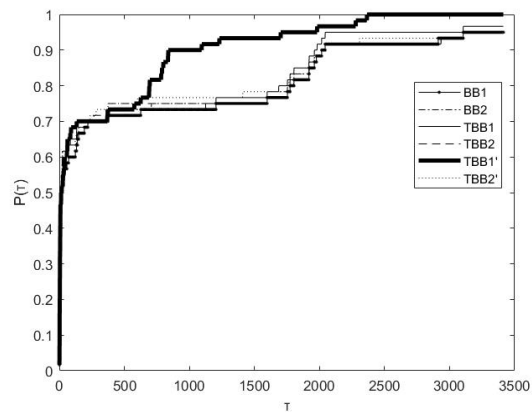

Figure 1 Performance evaluation of iteration times for 6 algorithms

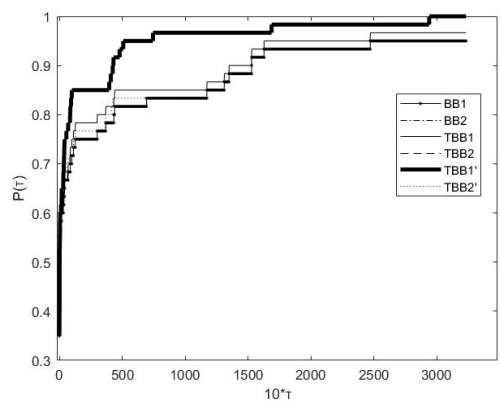

Figure 2 Performance evaluation of function value calculation times for 6 algorithms

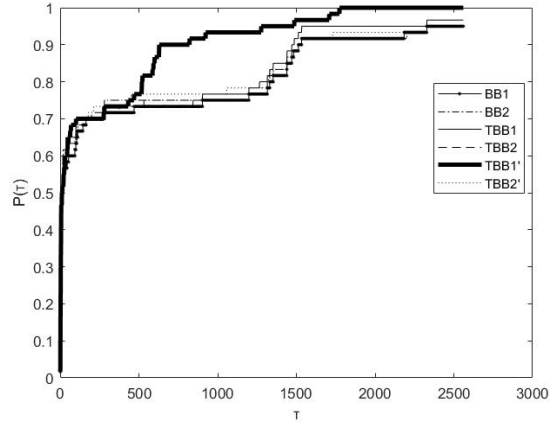

Figure 3 Performance evaluation of gradient calculation times for 6 algorithms

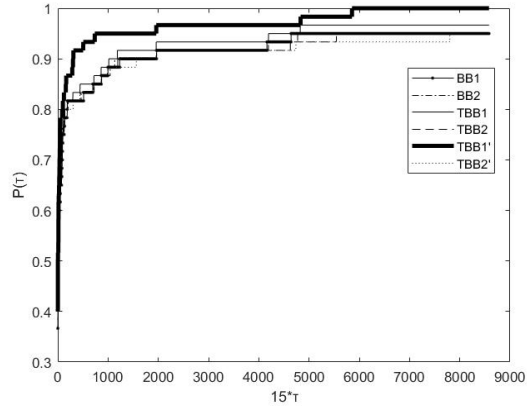

Figure 4 Performance evaluation of CPU time for 6 algorithms

In the second part, under the accuracy requirement $\|\nabla f(x_k)\| \leq 10^{-4}$, the results of Example 1 and Example 2 (N=10000, 100000, 1000000) are given in tables, where IT, T and $|f_{opt} - f(x^*)|$ represent the numbers of iteration, the CPU time and the function value error between the theoretical optimal solution and the real optimal solution, respectively. Also, the numerical results are given when $L_k = 0$.

**Example 8.1** Objective function $f = \sum_{i=1}^{N/2}[(x_{2i} - x_{2i-1}^2)^2 + (1 - x_{2i-1})^2]$

with initial point: $x_0 = (-1.2, 1, -1.2, 1, ..., -1.2, 1)^T$; optimal point: $x^* = (1, 1, ..., 1)^T$, optimal value $f^* = 0$. The computed results are shown in Tables 8.1-8.2.

**Table 8.1** Computed results based on Example 8.1

| method | N | IT | T/s | $|f_{opt} - f^*|$ |
|---|---|---|---|---|
| BB1 | 10000 | 24 | 0.036903 | 2.0368e-10 |
|  | 100000 | 25 | 0.174249 | 3.0565e-11 |
|  | 1000000 | 25 | 1.888609 | 3.0565e-10 |
| BB2 | 10000 | 28 | 0.018989 | 1.6407e-15 |
|  | 100000 | 28 | 0.124217 | 1.6407e-14 |
|  | 1000000 | 28 | 1.415083 | 1.6407e-13 |
| TBB1 | 10000 | 24 | 0.027272 | 1.9137e-09 |
|  | 100000 | 26 | 0.229340 | 5.5237e-09 |
|  | 1000000 | 27 | 1.723397 | 3.4134e-16 |
| TBB2 | 10000 | 20 | 0.020178 | 1.4653e-10 |
|  | 100000 | 22 | 0.168371 | 1.2974e-10 |
|  | 1000000 | 22 | 1.691725 | 1.2974e-09 |
| $TBB1'$ | 10000 | 24 | 0.025161 | 3.1672e-10 |
|  | 100000 | 24 | 0.162722 | 3.1672e-09 |
|  | 1000000 | 26 | 1.679969 | 1.0624e-08 |
| $TBB2'$ | 10000 | 20 | 0.019877 | 1.8264e-15 |
|  | 100000 | 20 | 0.139856 | 1.8264e-14 |
|  | 1000000 | 20 | 1.559532 | 1.8263e-13 |

**Table 8.2** Computed results based on Example 8.1 with $L_k = 0$

| method | N | IT | T/s | $|f_{opt} - f^*|$ |
|---|---|---|---|---|
| BB1 | 10000 | 45 | 0.069860 | 1.5664e-10 |
|  | 100000 | 45 | 0.254626 | 1.5664e-09 |
|  | 1000000 | 47 | 1.934844 | 8.7261e-09 |
| BB2 | 10000 | 44 | 0.087011 | 9.7986e-09 |
|  | 100000 | 45 | 0.332767 | 2.6498e-15 |
|  | 1000000 | 46 | 1.839784 | 2.6498e-14 |
| TBB1 | 10000 | 36 | 0.065912 | 8.2742e-09 |
|  | 100000 | 34 | 0.628646 | 1.0438e-08 |
|  | 1000000 | 35 | 1.543059 | 5.5892e-09 |
| TBB2 | 10000 | 40 | 0.022168 | 4.7195e-09 |
|  | 100000 | 41 | 0.174613 | 1.0243e-08 |
|  | 1000000 | 42 | 1.917521 | 4.7270e-17 |
| $TBB1'$ | 10000 | 38 | 0.025767 | 7.6468e-13 |
|  | 100000 | 38 | 0.194263 | 7.6468e-12 |
|  | 1000000 | 38 | 2.095670 | 7.6468e-11 |
| $TBB2'$ | 10000 | 30 | 0.029139 | 7.0063e-11 |
|  | 100000 | 31 | 0.163482 | 5.4858e-11 |
|  | 1000000 | 32 | 1.692133 | 8.9681e-11 |

From Tables 8.1-8.2 we see that as the numbers of variables become larger and larger, the numbers of iteration and CPU time increase slowly. This indicates that the six algorithms are efficient. The errors in the last columns show that the six algorithms are of global convergence.

**Example 8.2** Objective function

$$f(x) = \sum_{i=1}^{n/10}[(1-x_{10i-9})^2 + (1-x_{10i})^2 + \sum_{j=10i-9}^{10i-1}(x_j^2 - x_{j+1})^2].$$

with initial point: $x_0 = (-2,-2,\cdots,-2)^T$; optimal point: $x^* = (1,1,\cdots,1)^T$; optimal value: $f^* = 0$. The computed results are shown in Tables 8.3-8.4.

Table 8.3 Computed results based on Example 8.2

| method | N | IT | T/s | $\|f_{opt} - f^*\|$ |
|---|---|---|---|---|
| BB1 | 10000 | 37 | 0.019977 | 1.5102e-09 |
| | 100000 | 41 | 0.155796 | 1.8128e-09 |
| | 1000000 | 47 | 2.635774 | 5.8229e-11 |
| BB2 | 10000 | 39 | 0.019917 | 2.9408e-09 |
| | 100000 | 45 | 0.174077 | 9.7236e-10 |
| | 1000000 | 47 | 2.578728 | 1.2506e-11 |
| TBB1 | 10000 | 39 | 0.012575 | 1.5250e-10 |
| | 100000 | 39 | 0.147253 | 1.5250e-09 |
| | 1000000 | 42 | 2.349189 | 3.1141e-09 |
| TBB2 | 10000 | 38 | 0.013268 | 2.9423e-09 |
| | 100000 | 39 | 0.155520 | 1.4335e-09 |
| | 1000000 | 40 | 2.556326 | 2.2967e-10 |
| $TBB1'$ | 10000 | 40 | 0.016916 | 1.9239e-10 |
| | 100000 | 41 | 0.168539 | 1.9901e-10 |
| | 1000000 | 41 | 2.548831 | 1.9901e-09 |
| $TBB2'$ | 10000 | 38 | 0.019455 | 2.7199e-10 |
| | 100000 | 39 | 0.146647 | 7.4920e-10 |
| | 1000000 | 42 | 2.348636 | 1.1904e-09 |

Table 8.4 Computed results based on Example 8.2 with $L_k = 0$

| method | N | IT | T/s | $\|f_{opt} - f^*\|$ |
|---|---|---|---|---|
| BB1 | 10000 | 475 | 0.238484 | 1.4364e-09 |

|  | 100000 | 469 | 1.585360 | 9.1731e-11 |
|  | 1000000 | 473 | 24.247770 | 2.0115e-09 |
|  | 10000 | 270 | 0.093786 | 3.0030e-09 |
| BB2 | 100000 | 304 | 1.046501 | 2.9060e-09 |
|  | 1000000 | 276 | 14.753567 | 3.5457e-09 |
|  | 10000 | 483 | 0.191222 | 7.6750e-10 |
| TBB1 | 100000 | 523 | 1.946011 | 3.6674e-09 |
|  | 1000000 | 522 | 27.196727 | 1.4735e-09 |
|  | 10000 | 263 | 0.096177 | 9.4142e-10 |
| TBB2 | 100000 | 266 | 0.944474 | 8.9492e-10 |
|  | 1000000 | 279 | 15.401189 | 4.4790e-10 |
|  | 10000 | 515 | 0.200434 | 1.8464e-09 |
| $TBB1'$ | 100000 | 533 | 1.871263 | 2.0118e-09 |
|  | 1000000 | 534 | 27.949173 | 2.2781e-09 |
|  | 10000 | 302 | 0.122576 | 8.0590e-11 |
| $TBB2'$ | 100000 | 317 | 1.146826 | 2.2551e-09 |
|  | 1000000 | 297 | 15.904557 | 1.0218e-10 |

From the data in Tables 8.3-8.4 we also see that the numbers of iteration don't increase quickly but slowly as the problem dimension becomes larger. In fact, parameters in the four algorithms to solve different example problems can be selected differently, and the computed results will be improved. Here the same parameters are used to solve different optimization problems for simplicity.

By the computed results we see that Algorithms BB1, BB2，TBB1，TBB2，$TBB1'$ and $TBB2'$ are all globally convergent, that computed results by $TBB1'$ are better than two-point step size gradient algorithm. Upon search step size, relaxed non monotonic step size search is much better than monotonic step size search ($L_k = 0$). At approximation to optimal value $f^*$, the algorithms are of high accuracy. At the dimension of optimization problem, although the dimension increases quickly, the numbers of iteration and CPU time increase slowly. Thus the new algorithms proposed are suitable to solve large scale unconstrained optimization problems.

## 9. Conclusion

In this paper we proposed three-point step size gradient formula which extends the traditional BB formula，and then constructed relaxed generalized Armijo three-point step size gradient algorithms for solving unconstrained optimization problems. It is proved that the algorithms are globally convergent, linearly convergent and super linearly convergent in some cases. Moreover, it is shown that, when the objective function is pseudo-convex (quasi-convex) function, the new method has strong convergence results. It is demonstrated that the new proposed algorithms are simple in structure and robust, requiring less storage and information of first derivative. Numerical experiments show that the new algorithms are efficient to solve large scale unconstrained optimization problems.

**Acknowledgment**

This research is supported in part by Chinese Natural Foundation with grant number: 51974343. Corresponding author: Sun Qingying, Email: sunqingying01@163.com